\documentstyle[12pt, leqno]{article}

\begin{document}

\title{Symplectic reduction and a weighted multiplicity formula
 for  twisted Spin$^c$-Dirac operators}

\author{Youliang Tian\thanks{Partially supported by
a NYU research challenge fund grant.}
$\ $ and Weiping Zhang\thanks{Partially supported by  the NNSF,
SEC of China and the Qiu Shi Foundation.}}


\maketitle

\begin{abstract}
We extend our earlier work in [TZ1], where an analytic approach to the
 Guillemin-Sternberg conjecture [GS] was developed,
to cases where the Spin$^c$-complex under consideration is allowed
to be further twisted by certain natural
exterior power bundles. The main result is a {\it weighted}
quantization formula in the presence of commuting Hamiltonian actions.
The corresponding  Morse
type inequalities in  holomorphic situations are also established.
\end{abstract}

{\bf \S 0. Introduction}

$\ $

In a previous paper [TZ1], we have developed  a direct analytic approach to,
as well as certain extensions of,
the Guillemin-Sternberg geometric quantization conjecture [GS], which 
has been proved in various generalities in ([DGMW, G, GS, JK,
M1, M2, V1, V2]). In this paper, we generalize the results in
[TZ1] to  cases where the Spin$^c$-complex under consideration 
is allowed to be further twisted by certain natural
exterior power bundles. The main result is 
a {\it weighted} quantization formula for these twisted Spin$^c$-complexes
in the presence of commuting Hamiltonian actions. We also
establish the corresponding Morse type inequalities in the holomorphic
situation.

To begin with, let $(M,\omega)$ be a closed symplectic manifold 
admitting a Hamiltonian  action of a compact connected Lie group $G$
with Lie algebra ${\bf g}$.
Let $J$ be an almost complex structure on $TM$ so that $g^{TM}(u,v)
=\omega (u,Jv)$ defines a Riemannian metric on $TM$. After an integration
over $G$ if necessary, 
 we can and we will assume that $G$ preserves $J$ and $g^{TM}$.
Let $E$ be  a $G$-equivariant
Hermitian vector bundle with a $G$-equivariant 
Hermitian connection $\nabla^E$.

With these data in hand, for any integer $p\geq 0$, one can
construct canonically a formally self-adjoint  twisted
Spin$^c$-Dirac  operator acting 
on smooth sections of the twisted Spin$^c$-vector bundles:
$$D^{\wedge^{p,0}(T^*M)\otimes E}_+: 
\Omega ^{p,\rm even}(M,E)\rightarrow \Omega ^{p,\rm odd}(M,E).\eqno(0.1)$$
It gives rise to the finite dimensional virtual vector space
$$ Q(M,\wedge^{p,0}(T^*M)\otimes E) =\ker 
D^{\wedge^{p,0}(T^*M)\otimes E}_{+} - 
{\rm coker}\, D^{\wedge^{p,0}(T^*M)\otimes E}_{+}. \eqno(0.2)$$

Since $G$ preserves everything, one sees easily that 
$Q(M,\wedge^{p,0}(T^*M)\otimes E) $ is a  virtual representation
of $G$. Denote by $ Q(M,\wedge^{p,0}(T^*M)\otimes E)^G $ 
its $G$-trivial component.

Let ${\bf g}$ (and thus its dual ${\bf g}^*$ also) be equipped with
an Ad$\, G$-invariant metric. Let $h_i$, $1\leq i\leq \dim G$,
be an orthonormal base of ${\bf g}^*$.
Let $V_i$, $1\leq i\leq \dim G$, be the dual base of 
$\{ h_i\}_{1\leq i\leq \dim G}$.

Let $\mu:M\rightarrow {\bf g}^*$ be the moment map of the $G$ action on $M$.
Then it can be written as
$$\mu=\sum_{i=1}^{\dim G}\mu_ih_i, \eqno (0.3)$$
with each $\mu_i$ a real function on $M$.

Now for each $V\in {\bf g}$, set\footnote{When there is no confusion,
in this paper we will
use the same notation $V$ for the Killing vector field it generates on $M$.}
$$r_V^E=L_V^E-\nabla_V^E,\eqno (0.4)$$
where $L_V^E$ denotes the infinitesimal action of $V$ on $E$.

$\ $

{\bf Definition 0.1.} {\it We say $E$ is $\mu$-positive 
if the inequality
$$\sqrt{-1}\sum_{i=1}^{\dim G}\mu_i(x)r_{V_i}^E(x)> 0\eqno(0.5)$$
holds at every critical point 
$x\in M\backslash \mu^{-1}(0)$ of $|\mu|^2$, the norm square
of the moment map.}

$\ $

As a typical example, the  $G$-equivariant prequantum line bundle 
over $(M,\omega)$ verifying  the Kostant formula ([Ko], cf. [TZ1, (1.13)]),
when it exists,  is $\mu$-positive.

To state the main results of this paper, 
we now assume that $0\in {\bf g}^*$ is a
regular value of $\mu$ and, for simplicity,
that $G$ acts freely on $\mu^{-1}(0)$.
Then  one can construct the Marsden-Weinstein reduction
$(M_G,\omega_G)$, which is a smooth symplectic manifold with
$M_G=\mu^{-1}(0)/G$ and the 
symplectic form $\omega_G$ descended from $\omega$. The almost
complex structure $J$ also descends to an almost complex structure on
$TM_G$.
Furthermore, $E$ descends to a Hermitian vector bundle $E_G$ over $M_G$
with an induced Hermitian connection. Thus one can make the same
construction of the twisted Spin$^c$-Dirac operators as well
as the associated virtual vector spaces 
$Q(M_G,\wedge^{p,0}(T^*M_G)\otimes E_G)$.

For any integer $k,\ s\geq 0$, let $C_{s}^k $  be the binomial coefficient
given by
$$ C_{s}^k ={s(s-1)\cdots (s-k+1)\over
k!}.\eqno (0.6)$$

The main result of this paper, which is a generalization
of [TZ1, Theorem 4.1] in the abelian group action case, can be
stated as  follows.

$\ $

{\bf Theorem 0.2.} {\it If $G$ is abelian and $E$ is $\mu$-positive,
then the following identity holds for any
integer $p\geq 0$,
$$\dim Q(M,\wedge^{p,0}(T^*M)\otimes E)^G 
= \sum_{k=0 }^{p} 
C_{\dim G}^{p-k}\cdot \dim Q(M_G, 
\wedge^{k,0}(T^*M_G) \otimes E_G).\eqno (0.7)$$}

When $p=0$ and $E$ is the prequantum line bundle (when it exists)
over $(M,\omega)$, (0.7) is the abelian 
version of the Guillemin-Sternberg conjecture [GS]
proved by Guillemin [G] in a special case and by Meinrenken [M1] and 
Vergne [V1, 2] in general (see also [DGMW] and [JK]). Thus in some sense
one may view (0.7) as a kind of {\it weighted} quantization formula.

We will use the analytic approach developed in [TZ1] 
to prove Theorem 0.2. However, it
should be pointed out that Theorem 0.2  is not
a consequence of the  result in [TZ1, Theorem 4.1], which itself is a
generalization of the Guillemin-Sternberg conjecture [GS]. In
particular,  the strict inequality in (0.5) can not be relaxed to
include the equality as in [TZ1, Theorem 4.2], even
 when $\mu^{-1}(0)\neq \emptyset$. 
Furthermore, the abelian condition on $G$ is essential to both  the results 
as well as their proofs. 
A notable feature here is that we deal with the general case, where $G$
may possibly be of higher rank,  directly. 
That is, we do not first prove the result
for the $G=S^1$ case and then use the `reduction in stage' procedure
to get the full result.

Now  as in [TZ1, Theorem 0.4 and 4.8], we consider a holomorphic
refinement of Theorem 0.2. That is, we assume that $(M,\omega, J)$ is K\"ahler,
$G$ acts on $M$ holomorphically and $E$ is an equivariant
 holomorphic Hermitian
vector bundle over $M$ with the $G$ action on $E$ being holomorphic. If
for any integers $p,\ q\geq 0$,  denote by
$$h^{p,q}(E)^G=
  \dim H^{0,q}(M,\wedge^{p,0}(T^*M)\otimes E)^G,$$
$$h^{p,q}(E_G) =
 \dim H^{0,q}(M_G,\wedge^{p,0}(T^*M_G)\otimes E_G)\eqno (0.8)$$
the corresponding ($G$-invariant) twisted Hodge numbers, then we can state 
our refinement of (0.7) as follows.

$\ $

{\bf Theorem 0.3.} {\it If $G$ is abelian and $E$ is $\mu$-positive,
then the following inequality holds for any integers $p,\ q\geq 0$,
$$h^{p,q}(E)^G - h^{p,q-1}(E)^G+\cdots + (-1)^q h^{p,0}(E)^G$$
$$ \leq \sum_{k=0}^{p} C_{\dim G}^{p-k}
\left( h^{k,q}(E_G)- h^{k,q-1}(E_G)+\cdots
 +(-1)^q  h^{k,0}(E_G)\right).\eqno (0.9)$$
In particular, when $q=0$, one gets the following inequality for dimensions
of spaces of holomorphic sections,
$$h^{p,0}(E)^G \leq 
\sum_{k=0 }^{p} C_{\dim G}^{p-k}
\cdot h^{k,0}(E_G).\eqno (0.10)$$ }

Again, Theorem 0.3 is not a consequence of [TZ1, Theorem 4.8].

This paper is organized as follows. In Section 1, we construct the twisted
Spin$^c$-Dirac operators appearing in the context and introduce the 
corresponding deformations under Hamiltonian actions as in [TZ1].
We also prove a Bochner type formula for the deformed operators.
In Section 2, we extend the methods in [TZ1], which goes back to [BL],
to prove Theorems 0.2 and 0.3. The final Section 3
contains some immediate applications as well as further extensions of the
above main results.

$$\ $$

{\bf \S 1. Deformations of twisted Spin$^c$-Dirac operators and a Bochner
type formula}

$\ $

Following [TZ1], we construct
in this section  the twisted Spin$^c$-Dirac  operators  
and their deformations to be used in our proof of
Theorems 0.2 and 0.3. An important Bochner type formula
for the deformed operators will be proved.

This section is organized as follows. In a), we construct the above
mentioned Dirac  operators. In b), following [TZ1],
in the situations of Hamiltonian 
actions we introduce the deformations of the Dirac 
operators constructed in a). In c),  we prove a
Bochner type formula for the deformed operators. 

$\ $

{\bf a). Twisted Spin$^c$-Dirac operators on symplectic manifolds}

$\ $

Let $(M, \omega)$ be a closed symplectic manifold. Let $J$ be an
almost complex structure on $TM$ such that
$$g^{TM}(v,w)=\omega(v, Jw) \eqno(1.1)$$
defines a Riemannian metric on $TM$.
Let $TM_{\bf C}=TM\otimes {\bf C}$ denote the complexification of the
tangent bundle $TM$. Then one has the canonical (orthogonal) splittings
$$TM_{\bf C}=T^{(1,0)}M\oplus T^{(0,1)}M,$$
$$\wedge ^{*,*}(T^*M)=\bigoplus_{i,j=0}^{\dim _{\bf C}M}
\wedge ^{i,j}(T^*M),\eqno(1.2)$$
where
$$T^{(1,0)}M=\{ z\in TM_{\bf C}; Jz=\sqrt{-1} z\},$$
$$T^{(0,1)}M=\{z\in TM_{\bf C}; Jz=-\sqrt{-1} z\},$$
$$\wedge ^{i,j}(T^*M)=\wedge ^i(T^{(1,0)*}M)\otimes \wedge ^j
(T^{(0,1)*}M), \eqno(1.3)$$
and $\dim _{\bf C}M$ is the complex dimension of $M$.

For any $X\in TM$ whose complexification has the decomposition
$X=X_1+X_2\in T^{(1,0)}M\oplus T^{(0,1)}M$, let
${\bar X_1}^*\in T^{(0,1)*}M$ (resp. ${\bar X_2}^*\in T^{(1,0)*}M$)
 be the metric dual of $X_1$ (resp. $X_2$).  Set as in [BL, Sect. 5] that 
$$c(X)=\sqrt{2} {\bar X_1}^*\wedge \ -\sqrt{2}i_{X_2} .
\eqno (1.4)$$ 
Then $c(X)$ defines the
canonical Clifford action of $X$ on $\wedge^{0,*}(T^*M)$.
In particular, for any  $X,\, Y\in TM$, one has
$$c(X)c(Y) + c(Y)c(X) = -2 g^{TM}(X,Y).\eqno(1.5)$$

Let $\nabla^{TM}$ be the Levi-Civita connection of $g^{TM}$. Then  
$\nabla^{TM}$ together with the almost complex structure $J$ 
induce via projection a canonical Hermitian connection $\nabla^{T^{(1,0)}M}$
on $T^{(1,0)}M$. This in turn induces canonically, for any integer $p\geq 0$,
a  Hermitian connection $\nabla^{\wedge^{p,0}(T^*M)}$
on $\wedge^{p,0}(T^*M)$. On the other
hand, as was shown in [TZ1], $\nabla^{TM}$ 
 lifts canonically to a Hermitian connection $\nabla^{\wedge^{0,*}(T^*M)}$
on $\wedge^{0,*}(T^*M)$. Let $\nabla^{\wedge^{p,*}(T^*M)}$ be the
Hermitian connection on $\wedge^{p,*}(T^*M)$ obtained from the
tensor product of $\nabla^{\wedge^{p,0}(T^*M)}$ and 
$\nabla^{\wedge^{0,*}(T^*M)}$.

Now let $E$ be a Hermitian vector bundle over $M$ with a
Hermitian connection $\nabla^E$. 
Let  $\nabla^{ \wedge^{p,*}(T^*M)\otimes E }$ be the tensor product connection
of $\nabla^{ \wedge^{p,*}(T^*M) }$ and $\nabla^E$
on $\wedge^{p,*}(T^*M)\otimes E $.

Denote by $\Omega^{p,*}(M,E)$ the set of smooth sections of
$\wedge^{p,*}(T^*M)\otimes E$.

Let $e_1,\ldots, e_{\dim M}$ be an orthonormal base of $TM$. 

Following [TZ1], we now introduce the twisted Spin$^c$-Dirac  
operator we are interested in.

$\ $

{\bf Definition 1.1.} {\it The twisted Spin$^c$-Dirac operator
$D^{\wedge^{p,0}(T^*M)\otimes E}$  is defined by
$$ D^{\wedge^{p,0}(T^*M)\otimes E}  =\sum_{j=1}^{\dim M} c(e_j)
\nabla_{e_j}^{\wedge^{p,*}(T^*M)\otimes E}  :
\Omega^{p,*}(M,E) \rightarrow \Omega^{p,*}(M,E).\eqno(1.6)$$}

Clearly, $D^{\wedge^{p,0}(T^*M)\otimes E}$  is a formally self-adjoint
first order elliptic differential operator. Let
$D^{\wedge^{p,0}(T^*M)\otimes E}_+$   be the restriction of
$D^{\wedge^{p,0}(T^*M)\otimes E}$   on $\Omega^{p,\rm even}(M,E) $. Set
$$Q(M,\wedge^{p,0}(T^*M)\otimes E)
=\ker D^{\wedge^{p,0}(T^*M)\otimes E}_{+} - {\rm coker}\,
D^{\wedge^{p,0}(T^*M)\otimes E}_{+}. \eqno (1.7)$$

$\ $

{\bf b). Hamiltonian  actions and deformations of Dirac operators}

$\ $

Now suppose that $(M,\omega)$ admits a Hamiltonian action of
a compact connected Lie group $G$ with Lie algebra ${\bf g}$.
Let $\mu : M\rightarrow {\bf g}^*$
denote the corresponding moment map. After an integration 
over $G$ if necessary, 
we assume that $G$ preserves $J$ and thus also $g^{TM}$. We
also  assume that the $G$ action on $M$ lifts to a $G$ action on $E$
preserving the Hermitian metric as well as the Hermitian connection
$\nabla^E$ on $E$.

Let ${\bf g}$ (and thus ${\bf g}^*$ also) be
equipped with an Ad$\, G$-invariant metric.  Let ${\cal{H}} =|\mu |^2$ be the
norm square of the moment map $\mu$.
Then ${\cal{H}}$ is a $G$-invariant function on $M$. In particular, its
Hamiltonian vector field, denoted by $X^{{\cal{H}}}$, is $G$-invariant.
The following formula for $X^{{\cal{H}}}$ is clear,
$$X^{{\cal{H}}} =-J(d{\cal{H}}) ^*. \eqno(1.8)$$

Let $h_1,\ldots, h_{\dim G}$ be an orthonormal base of ${\bf g}^*$.
Then $\mu$ has the expression
$$\mu=\sum_{i=1}^{\dim G} \mu_ih_i,\eqno(1.9)$$
where each $\mu_i$ is a real valued function on $M$. Denote by $V_i$
the Killing vector field on $M$ induced by the dual of $h_i$.
One easily verifies that (cf. [TZ1, Sect. 1b)])
$$J(d\mu_i)^*=-V_i\eqno (1.10)$$ 
and
$$X^{\cal H}= -2J\sum_{i=1}^{\dim G} \mu_i(d\mu_i)^*
=2\sum_{i=1}^{\dim G}\mu_i V_i.\eqno (1.11)$$

We are now ready to introduce the crucial deformation following [TZ1].

$\ $

{\bf Definition 1.2.} {\it For any $T\in {\bf R}$,
let $ D^{\wedge^{p,0}(T^*M)\otimes E}_T $  be the operator defined by
$$D^{\wedge^{p,0}(T^*M)\otimes E}_T =D^{\wedge^{p,0}(T^*M)\otimes E}
 + {\sqrt{-1}T\over 2}c(X^{\cal H}) :\Omega^{p,*}(M,E)
\rightarrow \Omega^{p,*}(M,E).\eqno(1.12)$$}

Clearly, $D^{\wedge^{p,0}(T^*M)\otimes E}_T $ is a formally self-adjoint
first order elliptic differential operator.
Also, since $G$ preserves everything
and $X^{{\cal{H}}}$ is $G$-invariant, one sees that 
$ D^{\wedge^{p,0}(T^*M)\otimes E}_T $  is $G$-equivariant.
If we denote by $ D^{\wedge^{p,0}(T^*M)\otimes E}_{T,+} $  the
restriction of $ D^{\wedge^{p,0}(T^*M)\otimes E}_T $  on 
$\Omega^{p,\rm even}(M,E)$, then
$$Q_T(M,\wedge^{p,0}(T^*M)\otimes E)=
\ker D^{\wedge^{p,0}(T^*M)\otimes E}_{T,+} - {\rm coker}\,
D^{\wedge^{p,0}(T^*M)\otimes E}_{T,+} \eqno (1.13)$$ 
is a virtual $G$-representation. We use as usual
a superscript $G$ to denote its $G$-invariant component. 

Clearly, the following easy yet important identity holds 
for any $T\in {\bf R}$,
$$\dim Q_T(M,\wedge^{p,0}(T^*M)\otimes E)^G =
\dim Q(M, \wedge^{p,0}(T^*M) \otimes E)^G. \eqno (1.14)$$

$\ $

{\bf c). A Bochner type formula for $D^{\wedge^{p,0}(T^*M)\otimes E}_T $ }

$\ $

For any $V\in {\bf g}$, let $L_V$ denote the infinitesimal action induced
by $V$ on the corresponding vector bundles. We will in general omit
the superscripts of these bundles.
Let $ r_V^E $ be defined as in (0.4). 

For any $X,\ Y\in TM$ whose complexifications have the decomposition
$X=X_1+X_2\in T^{(1,0)}M\oplus T^{(0,1)}M$ and
$Y=Y_1+Y_2\in T^{(1,0)}M\oplus T^{(0,1)}M$ respectively, 
let $A(X,Y)$ be the endomorphism
of $\wedge^{p,0}(T^*M) $ defined by
$$A(X,Y)= {\bar X}^*_2\wedge i_{Y_1} .\eqno (1.15)$$

Let $e_1,\cdots,e_{\dim M}$ be an orthonormal base of $TM$.
Then one has the following  analogue of [TZ1, Lemma 1.5].

$\ $

{\bf Lemma 1.3.} {\it The following identity for operators acting on
$\Omega^{p,*}(M,E)$ holds,
$$L_V=\nabla_V+ r_V^E 
-{1\over 4} \sum_{j=1}^{\dim M} c(e_j)c(\nabla^{TM}_{e_j}V)
-{1\over 2} {\rm Tr}\left[\nabla^{T^{(1,0)}M}_. V\right] 
+\sum_{j=1}^{\dim M} A(e_j,\nabla^{TM}_{e_j}V).
\eqno (1.16)$$}

{\it Proof}. By proceeding as in [TZ1, Lemma 1.5], one sees easily that
one needs only to calculate $r_V^{\wedge^{p,0}(T^*M)}$.

Recall that $V$ acts on $TM$ by
$$L_V^{TM}X=\nabla^{TM}_VX- \nabla^{TM}_XV ,\ \ X\in \Gamma(TM),\eqno(1.17)$$
from which we have 
$$ r_V^{TM} (X)= -  \sum_{j=1}^{\dim M}  \langle \nabla^{TM}_XV ,e_j \rangle
e_j = \sum_{j=1}^{\dim M}  \langle \nabla^{TM}_{e_j }V ,X\rangle
e_j .\eqno (1.18)$$
From (1.18) one gets immediately that
$$ r_V^{T^*M} = \sum_{j=1}^{\dim M}  e_j^* \wedge i_{\nabla^{TM}_{e_j }V }.
\eqno (1.19)$$

By (1.19), (1.15) and the fact that the almost complex structure $J$ is
$G$-invariant, one deduces easily that for any integer $0\leq p\leq
\dim_{\bf C}M$,
$$r_V^{\wedge^{p,0}(T^*M)} = \sum_{j=1}^{\dim M}  A(e_j,\nabla^{TM}_{e_j}V).
\eqno (1.20)$$

(1.16) then follows from (1.20) and the arguments in
[TZ1, Lemma 1.5]. $\Box$

$\ $

We can now state the following analogue of [TZ1, Theorem 1.6].

$\ $

{\bf Theorem 1.4.}  {\it The following Bochner type formula holds,
$$D^{\wedge^{p,0}(T^*M)\otimes E,2}_T = D^{\wedge^{p,0}(T^*M)\otimes E,2}
 -2\sqrt{-1}T\sum_{i=1}^{\dim G} \mu_iL_{V_i}
+{T^2\over 4}\left|X^{{\cal{H}}} \right|^2 $$
$$+2\sqrt{-1}T\sum_{i=1}^{\dim G} \mu_ir_{V_i}^E 
+{T\over 2 }\sum_{i=1}^{\dim G} \left({\sqrt{-1} }c(JV_i)c(V_i)+|V_i|^2-
4\sqrt{-1}A(JV_i,V_i)\right)$$
$$+{\sqrt{-1}T\over 4}\sum_{j=1}^{\dim M}\left( c(e_j)
c(\nabla^{TM}_{e_j}X^{\cal H})+
4A(e_j,\nabla^{TM}_{e_j}X^{\cal H})\right) 
 - {\sqrt{-1}T\over 2} {\rm Tr}
\left[ \nabla^{T^{(1,0)}M}X^{\cal H}\right]
.\eqno(1.21)$$ }

{\it Proof}. As in [TZ1, (1.26) and (1.27)], one deduces from (1.12), (1.5)
that
$$D_T^2=D^2+{\sqrt{-1} T\over 2}\sum_{j=1}^{\dim M}c(e_j)c(\nabla_{e_j}
X^{\cal H})-{\sqrt{-1} T}
\nabla_{X^{\cal H}}+{T^2\over 4}\left|X^{{\cal{H}}} \right|^2 \eqno(1.22)$$ 
and, by using (1.10), (1.11) and Lemma 1.4, that
$$\nabla_{X^{\cal H}} =2\sum_{i=1}^{\dim G} \mu_i L_{V_i}
-2\sum_{i=1}^{\dim G} \mu_ir_{V_i}^E
+{1\over 4}\sum_{j=1}^{\dim M}c(e_j)c(\nabla_{e_j}X^{\cal H})$$
$$+{1\over 2} {\rm Tr}
\left[ \nabla^{T^{(1,0)}M}X^{\cal H}\right]
 -{1\over 2} \sum_{i=1}^{\dim G} c(JV_i)c(V_i)+{\sqrt{-1}\over 2}
\sum_{i=1}^{\dim G} |V_i|^2 $$
$$-\sum_{j=1}^{\dim M} A(e_j,\nabla_{e_j}X^{\cal H})
+ 2\sum_{i=1}^{\dim G}A(JV_i,V_i).\eqno (1.23)$$

(1.21) follows from (1.22) and (1.23).
$\Box$

$$\ $$

{\bf \S 2. Proof of the main theorems}

$\ $

In this section, we apply the methods and techniques in
[TZ1, Sects. 2-4], which are closely related to those in
[BL], to prove Theorems 0.2 and 0.3. As in [TZ1], the key technical point
is a pointwise estimate at each critical point 
$x\in M\backslash \mu^{-1}(0)$ of ${\cal H}=|\mu|^2$.

This section is organized as follows. In a), we prove the key pointwise
estimate mentioned above. In b), we prove Theorem 0.2 while Theorem 0.3
will be proved in c). 

$\ $

{\bf a). An estimate outside of $\mu^{-1}(0)$}

$\ $

Recall from Definition 0.1 that $E$ is said to be $\mu$-positive if (0.5)
holds at every critical point $x\in M\backslash \mu^{-1}(0)$ of 
${\cal H}=|\mu|^2$.

The main result of this subsection, which is an analogue of [TZ1, Theorem 2.1],
can be stated as follows.

$\ $

{\bf Theorem 2.1.} {\it If $G$ is abelian and 
$E$ is $\mu$-positive, then
for any open neighborhood $U$ of $\mu^{-1}(0)$,
there exist constants $C>0,\ b>0$ such that for any $T\geq 1$
and any $G$-invariant section $s\in \Omega^{p,*}(M,E)$  
with ${\rm Supp}\; s\subset M\backslash U$,
one has the following estimate of Sobolev norms,}
$$\| D^{\wedge^{p,0}(T^*M) \otimes E}_T s\|^2_{0}
\geq C(\| s\|^2_{1} + (T-b)\| s\|^2_{0}).\eqno(2.1)$$

$Proof$. By examining the arguments in [TZ1, Sect. 2], one sees that 
in order to prove Theorem 2.1, one needs only to prove an
analogue of [TZ1, Lemma 2.3] in our context.

Thus let $x\in M\backslash \mu^{-1}(0)$ be a critical point of ${\cal H}$.
Let $e_1, \cdots ,e_{\dim M}$ be an orthonormal base of $TM$ near $x$. Let 
$(y_1, \cdots, y_{\dim M})$ be the normal coordinate system with respect to
$\{ e_j\} _{j=1}^{\dim M}$ near $x$.
Clearly, one can choose  $e_1, \cdots ,e_{\dim M}$ 
so that ${\cal{H}}$ has the following expression near $x$,
$${\cal{H}} (y)={\cal{H}} (x)+\sum_{j=1}^{\dim M} 
a_jy_j^2 +O(|y|^3),\eqno(2.2)$$
where the $a_j$'s may possibly be zero.

We can now state our analogue of [TZ1, Lemma 2.3] as follows.

$\ $

{\bf Lemma 2.2.} {\it If $G$ is abelian,
 then the following inequality holds at any 
critical point $x\in M\backslash \mu^{-1}(0)$ of ${\cal H}$,
$$ {\sqrt{-1}\over 4}\sum_{j=1}^{\dim M}\left( c(e_j)
c(\nabla^{TM}_{e_j}X^{\cal H})+
4A(e_j,\nabla^{TM}_{e_j}X^{\cal H})\right) 
- {\sqrt{-1}\over 2} {\rm Tr}
\left[ \nabla^{T^{(1,0)}M}X^{\cal H}\right] $$
 $$ + {1\over 2 }\sum_{i=1}^{\dim G} \left({\sqrt{-1} }c(JV_i)c(V_i)+|V_i|^2-
4\sqrt{-1}A(JV_i,V_i)\right) 
\geq -\sum_{j=1}^{\dim M}|a_j|. \eqno (2.3)$$}

{\it Proof}. Since $x$ is a critical point of ${\cal H}$, by a result of
Kirwan [K, Prop. 3.12], $x$ is a fixed point
for the action of the subtorus generated by $\mu(x)\neq 0$. 
Without loss of generality, we 
assume that $h_1, \cdots, h_{\dim G}$
has been chosen so that the duals of
$h_1, \cdots, h_r$ generate 
this subtorus denoted by $G_1$. Let $G_2$ be the subtorus
generated by the duals of $h_{r+1}, \cdots, h_{\dim G}$. 
Thus  the original torus has the splitting
$$G=G_1\times G_2.\eqno(2.4)$$
Clearly, one has that
$$ \mu_i(x)=0,\ \ r+1\leq i\leq \dim G.\eqno(2.5)$$

Denote by $F_x\subset M$ the connected component containing $x$
of the fixed point set of the $G_1$-action. 
Then  $F_x$ is a totally geodesic submanifold of
$M$ and $J$ preserves the tangent bundle $TF_x$. Denote $k=\dim F_x$. 

Since the $G_2$-action commutes with the $G_1$-action, $G_2$ acts on
$F_x$. 

To summarize, one has 

$\ $

{\bf Lemma 2.3.} {\it a). If $1\leq i\leq r$, then
$V_i(x)=0$ and $ \mu_i|_{F_x}$ is  constant;\\
b). if $r+1\leq i\leq \dim G$, then
$\mu_i(x)=0$ and  $(d\mu_i)|_{F_x}\in \Gamma (T^*F_x)$.}

{\it Proof}. Since $G_2$ acts on $F_x$, for any $r+1\leq i\leq \dim G$,
$V_i|_{F_x}\in \Gamma (TF_x)$. Thus by (1.10), 
$(d\mu_i)^*=JV_i\in \Gamma (TF_x)$. 
The other parts of the lemma are clear. $\Box$

$\ $

Without loss of generality, one may choose 
$e_1,\cdots, e_{\dim M}$ near $x$ so that
$ e_1, \cdots,e_k$ 
is an orthonormal base of  $TF_x$. Let $\{ y_j\}_{1\leq j\leq \dim M}$ 
be the corresponding normal coordinates near $x$. Then from Lemma 2.3
one deduces that, near $x$,  this orthonormal base can be further arranged 
so that
$$\sum_{i=1}^r |\mu_i(y)|^2 = |\mu(x)|^2
+\sum_{j=k+1}^{\dim M} a_j y_j^2 + O(|y|^3)\eqno (2.6)$$
and
$$\sum_{i=r+1}^{\dim G} |\mu_i(y)|^2 = 
\sum_{j=1}^k  a_j y_j^2+ O(|y|^3).\eqno(2.7)$$
(2.6) and (2.7) together provide a splitting 
of ${\cal H}$ near $x$ according to the splitting  (2.4).

From (1.10) and (2.6) one deduces, at $x$, that
$$\sum_{i=1}^{r}\nabla^{TM}_{e_j}( \mu_i V_i)=\Big\{ \matrix{ -a_j Je_j,&
{\rm for}\ k+1\leq j\leq \dim M, \cr 0,\;& {\rm for }\ 
1\leq j\leq k.\cr}\eqno(2.8)$$
From (2.8), one deduces, at $x$, that
$$ {\sqrt{-1}\over 4} \sum_{j=1}^{\dim M}\sum_{i=1}^r \left( c(e_j)
c( 2\nabla^{TM}_{e_j} (\mu_iV_i) )
+4A(e_j,2\nabla^{TM}_{e_j}(\mu_iV_i))\right) $$
$$-{\sqrt{-1}\over 2}\sum_{i=1}^r {\rm Tr}
\left[ 2\nabla^{T^{(1,0)}M}(\mu_iV_i)  \right] $$
$$=-\sum_{j=k+1}^{\dim M} a_j \left( {\sqrt{-1}\over 2} c(e_j)c(Je_j)
+ 2\sqrt{-1}A(e_j,Je_j) +{1\over 2} \right)
\geq -\sum_{j=k+1}^{\dim M} |a_j|,\eqno (2.9)$$
where the last inequality follows from the obvious inequalities that
$$\left| c(e_j)c(Je_j) \right| \leq 1\eqno (2.10)$$
and
$$\left| 2\sqrt{-1}A(e_j,Je_j) + {1\over 2} \right| =
\left|  {1\over 2} -2\overline{e_j^{0,1}}^*\wedge
i_{ e_j^{1,0} }\right| \leq {1\over 2} ,\eqno (2.11)$$
with $e_j^{1,0} \in T^{(1,0)}M$ (resp. $e_j^{0,1} \in T^{(0,1)}M$)
the (1,0) (resp. (0,1)) component of the complexification of $e_j$. 

On the other hand, by Lemma 2.3,
for each $r+1\leq i\leq \dim G$,  $\mu_i$ can
be written, near $x$, as
$$\mu_i(y)=\sum_{j=1}^{k}b_{ij}y_j+O(|y|^2).\eqno(2.12)$$

By (2.7) and (2.12), one deduces that
$$\sum_{i=r+1}^{\dim G}\sum_{j=1}^{k}b_{ij}^2=\sum_{j=1}^{k}
a_j,\eqno (2.13)$$
which, together with (1.10), imply 
$$\sum_{i=r+1}^{\dim G}|V_i(x)|^2 =\sum_{j=1}^k a_j.\eqno(2.14)$$

From Lemma 2.3, (2.14), (1.10) and (2.10), 
one deduces, at $x$, that
$$ {\sqrt{-1}\over 2} \sum_{j=1}^{\dim M}\sum_{i=r+1}^{\dim G} \left( c(e_j)
c( \nabla^{TM}_{e_j} (\mu_iV_i) )
+4A(e_j,\nabla^{TM}_{e_j}(\mu_iV_i))\right) $$
$$+ \sum_{i=r+1}^{\dim G}  \left(-{\sqrt{-1}}{\rm Tr}
\left[ \nabla^{T^{(1,0)}M}(\mu_iV_i)  \right] +
{\sqrt{-1}c(JV_i)c(V_i)+|V_i|^2\over 2 } \right)$$
$$-2\sqrt{-1} \sum_{i=r+1}^{\dim G}  A(JV_i,V_i) 
=\sum_{i=r+1}^{\dim G} \sqrt{-1} c(JV_i)c(V_i) 
\geq -\sum_{j=1}^{k}| a_j|.\eqno (2.15)$$

From (1.11), (2.9), (2.15) and Lemma 2.3, one gets (2.3). The proof
of Lemma 2.2 is completed. $\Box$

$\ $

Since $E$ is $\mu$-positive,
 using Theorem 1.4, (0.5) and Lemma 2.2, one can proceed in  the
same way as in [TZ1, Sect. 2], with almost
no changes, to prove Theorem 2.1. That is, we first
prove pointwise estimates analogous to [TZ1, Prop. 2.2]
 around each point outside of $\mu^{-1}(0)$, in using
Lemma 2.2 when dealing with critical points of ${\cal H}$, and then glue them together
to get the global estimate (2.1). The essential point in this last
gluing step is again as in [TZ1] that each $L_{V_i}$, $1\leq i\leq \dim G$,
vanishes when acting on $G$-invariant sections.
We leave the easy details to the interested reader. $\Box$

$\ $

{\bf Remark 2.4.} A notable
difference between Lemma 2.2 and [TZ1, Lemma 2.3] is that
even with some of the $a_j$'s  being negative, one does not have in general
a strict inequality in (2.3). This means that the $\mu$-positivity of
$E$ is necessary for Theorem 2.1 (compare with [TZ1, Remark 2.4]).

$\ $

{\bf b). A proof of Theorem 0.2}

$\ $

If $F$ is a $G$-equivariant Hermitian vector bundle over $M$,
we denote by $F_G$ its induced bundle on $M_G$ (cf. [TZ1, Sect. 4a)]).

As in [TZ1], Theorem 2.1 allows us to localize our problem to sufficiently
small neighborhoods of $\mu^{-1}(0)$. While near $\mu^{-1}(0)$, 
we can directly apply the analysis and results in [TZ1, Sects. 3 and 4a)],
which are closely related to [BL, Sects. 8, 9],
to the $G$-invariant restriction of $D^{ \wedge^{p,0}(T^*M) \otimes E }_T$.

Combining the above  arguments, one deduces  using  (1.14) the following
analogue of [TZ1, (4.3)] (of course for different twisted bundles and
conditions),
$$ \dim Q(M,\wedge^{p,0}(T^*M)\otimes E)^G = 
\dim Q(M_G, (\wedge^{p,0}(T^*M) \otimes E)_G ).
\eqno (2.16)$$

Now since $G$ is abelian, the normal bundle to $\mu^{-1}(0)$ is
equivariantly trivial. From this fact one deduces directly the 
{\it weighted} splitting
$$(\wedge^{p,0}(T^*M) )_G =\bigoplus_{k=0}^p\, C_{\dim G}^{p-k}
\cdot \wedge^{k,0}(T^*M_G),\eqno (2.17)$$
where the numbers $C_{\dim G}^{p-k}$ have been defined in (0.6).

Theorem 0.2 then follows from (2.16) and (1.17). $\Box$

$\ $

{\bf Remark 2.5.} In view of Remark 2.4, the 
$\mu$-positivity condition (0.5) can not be
weakened in general to include the equality as in [TZ1, Theorem 4.2],
even when $\mu^{-1}(0)$ is nonempty.

$\ $

{\bf Remark 2.6.} Though its proof is of the same method, Theorem 0.2
can not be deduced from results in [TZ1] without
imposing further conditions. The point is that if one wants to apply directly
the results in [TZ1] to our situation, one needs the condition that
at every critical point $x\in M\backslash \mu^{-1}(0)$ of
${\cal H}$,
$$\sqrt{-1}\sum_{i=1}^{\dim G}\mu_ir_{V_i}^{\wedge^{p,0}(T^*M) \otimes E}
\geq 0,\eqno (2.18)$$
which clearly is not a consequence of (0.5).

$\ $

{\bf c). A proof of Theorem 0.3}

$\ $

We now assume that $(M,\omega, J)$ is K\"ahler and $G$ acts on $M$ 
holomorphically. Furthermore, we assume that $E$ is a 
$G$-equivariant holomorphic 
Hermitian vector bundle over $M$ on which $G$ acts holomorphically and
that $\nabla^E$ is the unique holomorphic Hermitian connection.

The key observation is, similarly as in [TZ1, Remark 1.4], 
that for any $T\in {\bf R}$ we have in this situation
$$D^{ \wedge^{p,0}(T^*M) \otimes E }_T=\sqrt{2}
( e^{-T{\cal H}/2} \bar{\partial}^{\wedge^{p,0}(T^*M) \otimes E }
e^{T{\cal H}/2} +e^{T{\cal H}/2} 
(\bar{\partial}^{\wedge^{p,0}(T^*M) \otimes E })^*e^{-T{\cal H}/2} ).
\eqno (2.19)$$
Furthermore, one has clearly  an analogue of [TZ1, (3.54)].
Thus by the same reason as in [TZ1, Sect. 4d)], all the arguments
before this subsection preserve the ${\bf Z}$-grading nature of the 
twisted Dolbeault complex on $M$ with coefficient bundle 
$\wedge^{p,0}(T^*M)\otimes E$, and this leads to the following
holomorphic refinement of (2.16) which holds for any integer $q\geq 0$,
$$h^{p,q}(E)^G - h^{p,q-1}(E)^G+\cdots + (-1)^q h^{p,0}(E)^G
 \leq h^{0,q}( (\wedge^{p,0}(T^*M) \otimes E)_G )$$
$$- h^{0,q-1}((\wedge^{p,0}(T^*M) \otimes E)_G )
+\cdots  +(-1)^q  h^{0,0}((\wedge^{p,0}(T^*M) \otimes E)_G ).\eqno (2.20)$$

On the other hand, one also verifies directly that the splitting (2.17)
is holomorphic in this situation. 

(0.9) is then a consequence of
(2.17) and (2.20).  $\Box$

$$\ $$

{\bf \S 3. Applications and further extensions}

$\ $

In this section, we discuss some immediate applications and
possible extensions of Theorems 0.2, 0.3 as well as the methods and 
techniques involved in their proofs.

This section is organized as follows. In a), we apply Theorem 0.2 to get
a vanishing multiplicity result for twisted de Rham-Hodge 
operators and a weighted multiplicity formula for twisted Signature 
operators. In b), we prove a negative analogue of Theorem 0.2, that is, 
we prove a weighted multiplicity formula in the case that `$>$' is
replaced by `$<$' in (0.5). We also show that the strict inequalities are 
necessary. In c), we discuss briefly the case where $0\in {\bf g}^*$ is
not a regular value of the moment map $\mu$.     
Finally, we discuss in d) the applications to the typical example
where $E$ is the prequantum line bundle over $(M,\omega)$.

$\ $

{\bf a). Applications to twisted de Rham-Hodge and Signature operators}

$\ $

Set 
$$ Q_{\rm dR}(M,E) = \bigoplus_{p=\rm even} Q(M, \wedge^{p,0}(T^*M) \otimes E)
- \bigoplus_{p=\rm odd} Q(M, \wedge^{p,0}(T^*M) \otimes E)\eqno (3.1)$$
and 
$$ Q_{\rm Sig}(M,E) = \bigoplus_{p=0}^{\dim_{\bf C}M} 
Q(M, \wedge^{p,0}(T^*M) \otimes E). \eqno (3.2)$$
One verifies easily that $Q_{\rm dR}(M,E) $ (resp. $Q_{\rm Sig}(M,E) $)
is exactly the virtual vector space associated to the twisted (by $E$)
de Rham-Hodge (resp. Signature) operator on $M$. The following result
gives the corresponding multiplicity formulas for these operators.

$\ $

{\bf Theorem 3.1.} {\it Under the same assumptions  as in Theorem 0.2, 
the following identities hold,
$$\dim Q_{\rm dR}(M,E)^G =0,\eqno (3.3)$$
$$\dim Q_{\rm Sig}(M,E)^G =2^{\dim G} \dim Q_{\rm Sig}(M_G,E_G).
\eqno (3.4)$$}

{\it Proof}. Theorem 3.1 follows easily from Theorem 0.2 and the definitions
(3.1), (3.2) with some elementary computation. $\Box$

$\ $

{\bf Remark 3.2.} The assumption that
 $0\in {\bf g}^*$ is a regular value of $\mu$ 
is essential, particularly for the vanishing property  (3.3). 
This will be discussed further in c) and d).

$\ $

{\bf b). Weighted multiplicity formula for $\mu$-negative bundles}

$\ $

A $G$-equivariant Hermitian vector bundle $E$ over $(M,\omega)$ with 
$G$-invariant Hermitian
connection $\nabla^E$ is said to be $\mu$-negative if at every critical
point $x\in M\backslash \mu^{-1}(0)$ of ${\cal H}=|\mu|^2$, one has
$$ \sqrt{-1}\sum_{i=1}^{\dim G}\mu_i(x)r_{V_i}^E(x) < 0,\eqno(3.5)$$ 
instead of (0.5).

For any $\mu$-negative bundle $E$, one can introduce the same deformation
of twisted Spin$^c$-Dirac operators as in Definition 1.2, but take
$T\rightarrow -\infty$, instead of $+\infty$, to prove the following result.

$\ $

{\bf Theorem 3.3.} {\it If $G$ is abelian and $E$ is $\mu$-negative, then
the following identity holds for any integer $p\geq 0$,
$$\dim Q(M,\wedge^{p,0}(T^*M)\otimes E)^G $$
$$= (-1)^{\dim G} \sum_{k=0 }^{p} C_{\dim G}^{p-k}\cdot \dim Q(M_G, 
\wedge^{k,0}(T^*M_G) \otimes E_G).\eqno (3.6)$$}

{\it Proof}. On can proceed similarly as in Section 2 to prove Theorem 3.2.
The key point to note is that now the analogue of Lemma 2.2 should take
the following form at any critical point 
$x\in M\backslash \mu^{-1}(0)$ of ${\cal H}=|\mu|^2$, 
$$ {\sqrt{-1}\over 4}\sum_{j=1}^{\dim M}\left( c(e_j)
c(\nabla^{TM}_{e_j}X^{\cal H})+
4A(e_j,\nabla^{TM}_{e_j}X^{\cal H})\right) 
- {\sqrt{-1}\over 2} {\rm Tr}
\left[ \nabla^{T^{(1,0)}M}X^{\cal H}\right] $$
 $$ + {1\over 2 }\sum_{i=1}^{\dim G} \left({\sqrt{-1} }c(JV_i)c(V_i)+|V_i|^2-
4\sqrt{-1}A(JV_i,V_i)\right) 
\leq \sum_{j=1}^{\dim M}|a_j|. \eqno (3.7)$$

We leave the details to the interested reader. $\Box$

$\ $

As immediate applications, one  gets analogues of Theorem 3.1 for
$\mu$-negative bundles. One also gets Morse type inequalities as
refinements of (3.6) in  holomorphic situations.

$\ $

{\bf Remark 3.4.} We use a simple example to illustrate that 
even with $\mu^{-1}(0)\neq \emptyset$, one can not weaken the 
$\mu$-positive (resp. $\mu$-negative) assumption in (0.5) (resp. (3.5))
by an equality to get 
a weighted multiplicity formula similar to what happens in
[TZ1, Theorem 4.2]: Taking $M$ to be ${\bf P}^1$ with its standard
$S^1$ action and $E={\bf C}$,  one verifies directly that
$\dim Q_{\rm dR}(M,{\bf C}) ^G=\dim Q_{\rm dR}(M,{\bf C}) =2\neq 0$.

$\ $

{\bf Remark 3.5.} It is also clear that
the abelian condition on $G$ is essential for our argument. In fact,
an example due to Vergne (cf. [JK, pp.686]) shows that the abelian condition
on $G$ is necessary for Theorem 3.3 in the case where
$p=0$ and $E$ is the dual
of the prequantum line bundle over $(M,\omega)$ (when the latter exists).

$\ $

{\bf c). The case where $0\in {\bf g}^*$ is a singular value of $\mu$}

$\ $

We now make a brief discussion on the possible generalizations of our
main results to the case where $0\in {\bf g}^*$ is a singular value of $\mu$.

When $p=0$ and $E$ is the prequantum line bundle over
$(M,\omega)$ (when it exists), 
the  quantization formula for singular reductions has been established
 by Meinrenken-Sjamaar [MS] (see also [TZ3] for an 
 analytic treatment in certain  situations).
One notable feature in this case is the phenomenon that `the singular value $0$
is removable' in the perturbative singular
quantization formula [MS, Theorems 2.5, 2.9] (see also [TZ3, Theorem 0.1]). 
However, as we will  see, the situation for
the case where $p$ is nonzero is rather different.

In the simplest case where  $G$ is the circle, a fairly general singular
localization formula, which can indeed be applied to our situation,
has been proved in [TZ2, Theorem 6.7].  To be more precise, let $V$ be
the Killing vector field on $M$ generated by the unit base of ${\bf g}$
and let $F_0= \mu^{-1}(0)\cap \{ x\in M; V(x)=0\}$ be the subset of
the fixed points of the $G$ action contained in $\mu^{-1}(0)$. Then  using
the same notation as in [TZ2, Sect. 6c) and Appendix], one has the following
result, which can be proved by a   combination of the arguments in
Section 2 with those in the proof of [TZ2, Theorem 6.7].

$\ $

{\bf Theorem 3.6.} {\it If $G=S^1$, $E$ is $\mu$-positive and
$0\in {\bf g}^*$ is a singular value of $\mu$, then
$$\dim Q(M, \wedge^{p,0}(T^*M)\otimes E)^G $$
$$= \sum_{k=0 }^{p} C_{\dim G}^{p-k}\cdot 
\dim Q(M_{G,0^-}, \wedge^{k,0}(T^*M_{G,0^-})\otimes E_{G,0^-}) +
{\rm ind}\, D_{F_0,+}^{\wedge^{p,0}(T^*M)\otimes E}(V) $$
$$= \sum_{k=0 }^{p} C_{\dim G}^{p-k}\cdot 
\dim Q(M_{G,0^+}, \wedge^{k,0}(T^*M_{G,0^+})\otimes E_{G,0^+}) +
{\rm ind}\, D_{F_0,+}^{\wedge^{p,0}(T^*M)\otimes E}(-V) .\eqno (3.8)$$}

$\ $

{\bf Corollary 3.7.} {\it Under the same conditions as in Theorem 3.6,
one has
$$ \dim Q_{\rm dR}(M,E)^G =
{\rm ind}\, D_{F_0,+}^{\wedge^{{\rm even},0}(T^*M)\otimes E}(V)
- {\rm ind}\, D_{F_0,+}^{\wedge^{{\rm odd},0}(T^*M)\otimes E}(V) $$
$$={\rm ind}\, D_{F_0,+}^{\wedge^{{\rm even},0}(T^*M)\otimes E}(-V)
- {\rm ind}\, D_{F_0,+}^{\wedge^{{\rm odd},0}(T^*M)\otimes E}(-V).
\eqno (3.9)$$}

{\bf Remark 3.8.} In the next subsection, we will show that even when $E$
is the prequantum line bundle over $(M,\omega)$, the contributions of 
$F_0$ to
the right hand sides of (3.8), (3.9) may well be nonzero. This explains the
essential difference  between  the  $p=0$ and $p\neq 0$ 
cases as mentioned above.

$\ $

{\bf Remark 3.9.} If $G$ is of higher rank, 
one can apply Theorems 3.6, 3.7 inductively
to get localization formulas for $\dim Q(M, \wedge^{p,0}(T^*M)\otimes E)^G $
and $\dim Q_{\rm dR}(M,E)^G $ respectively.

$\ $

{\bf d). The case where $E$ is the prequantum line bundle over
$(M,\omega)$}

$\ $

In this section, we assume $E$ is the prequantum line bundle $L$ over
$(M,\omega)$, of which we assume the existence. Then $L$ is a $G$-equivariant
Hermitian vector bundle over $M$ with  $G$-invariant Hermitian
connection $\nabla^L$ such that 
$${\sqrt{-1} \over 2\pi} \left( \nabla^L\right) ^2 =\omega .\eqno (3.10)$$
Furthermore, we assume that $\mu$ verifies the  
Kostant formula  [Ko] (cf. [TZ1, (1.13)])
$${ L}_V s=\nabla^L_V s-2\pi\sqrt{-1}\langle \mu,V\rangle s, 
\ s\in \Gamma (L),\ V\in {\bf g}. \eqno(3.11)$$

The following result is clear.

$\ $

{\bf Proposition 3.10.} {\it $L$ is $\mu$-positive.}

$\ $

One of the  novelties for the prequantum line bundle is that there is a standard
shifting trick
 to reduce the computation of dimensions of nontrivial components of the 
$G$-representation $Q(M, \wedge^{p,0}(T^*M)\otimes L)$ to those of 
trivial components of $G$-representations of the form
 $Q(M, \wedge^{p,0}(T^*M)\otimes E)$ with $E$ to be `$\mu$-positive'.

To be more precise, for any $\xi_i\in {\bf Z}$, $1\leq i\leq \dim G$,
set 
$$\xi=\xi_1h_1+\cdots \xi_{\dim G}h_{\dim G}\in {\bf g}^*.\eqno (3.12)$$
Let $ {\bf C}_{\xi} $ be the $G$-equivariant trivial complex line bundle
over $M$ with $G$-invariant Hermitian connection $\nabla^{{\bf C}_{\xi} } $,
on which $G$ acts by
$${ L}_V s= \nabla^{{\bf C}_{\xi} } _V s+2\pi\sqrt{-1}\langle \xi,V\rangle s, 
\ s\in \Gamma (L),\ V\in {\bf g}. \eqno(3.13)$$
The existence of $ {\bf C}_{\xi} $ is clear.

Set
$$ L_{\xi} =L\otimes {\bf C}_{\xi} .\eqno (3.14)$$
Then $L_{\xi} $ is canonically a $G$-equivariant Hermitian vector bundle
over $M$ with the tensor product connection $\nabla^{L_{\xi} }$.
Furthermore, $G$ acts on $L_{\xi} $ through the formula
$${ L}_V s= \nabla^{L_{\xi} } _V s- 2\pi\sqrt{-1} \langle \mu- \xi,V\rangle s, 
\ s\in \Gamma (L_{\xi}),\ V\in {\bf g}. \eqno(3.15)$$

Since $G$ is abelian, $\mu_{\xi}:=\mu-\xi$ may also be regarded as a 
moment map for the Hamiltonian action of $G$ on $(M,\omega)$. One then has the
following extension of Proposition 3.10, which can be verified directly.

$\ $

{\bf Proposition 3.11.} {\it $L_{\xi}$ is $\mu_{\xi}$-positive.}

$\ $

Now let $ Q(M, \wedge^{p,0}(T^*M)\otimes L)_{\xi} $ denote 
the $\xi$-eigen-component of the $G$-representation 
$Q(M, \wedge^{p,0}(T^*M)\otimes L)$.
That is, $L_V$ acts on $Q(M, \wedge^{p,0}(T^*M)\otimes L)_{\xi}$ by
multiplication by $2\pi\sqrt{-1} \langle \xi,V\rangle $. Then one
verifies directly that
$$ \dim Q(M, \wedge^{p,0}(T^*M)\otimes L) _{\xi} =
\dim Q(M, \wedge^{p,0}(T^*M)\otimes L_{\xi})^G .\eqno (3.16)$$ 
This is the shifting trick   mentioned above.

Now when $\xi$ is a regular value of $\mu$,  in view of
Proposition 3.11 one can apply Theorem 0.2
to get a weighted multiplicity formula calculating 
$\dim Q(M, \wedge^{p,0}(T^*M)\otimes L_{\xi})^G $ and thus the
$\xi$-multiplicity of $ Q(M, \wedge^{p,0}(T^*M)\otimes L) $ in terms
of quantities on the $\mu_{\xi}$-symplectic reduction
$M_{G,\xi}=\mu^{-1}(\xi)/G$. While when $\xi$ is not 
a regular value of $\mu$,  one may first apply Theorem 3.6 and then an
induction procedure to calculate 
$\dim Q(M, \wedge^{p,0}(T^*M)\otimes L) _{\xi} $ via (3.16).

By summing over all these $\xi$'s, one has clearly that
$$ \dim Q(M, \wedge^{p,0}(T^*M)\otimes L) =\sum_{\xi}
\dim Q(M, \wedge^{p,0}(T^*M)\otimes L) _{\xi} .\eqno (3.17)$$
In particular, in the case where $G=S^1$,
if for any such $\xi$ denote by $F_{\xi}=\mu^{-1}(\xi)\cap F$
where $F$ is the fixed point set of the $G$ action on $M$,
then one can combine the above reasoning with Corollary 3.7 to get
$$ \dim Q_{\rm dR}(M,L) =\sum_{\xi} \left(
{\rm ind}\, D_{F_{\xi},+}^{\wedge^{{\rm even},0}(T^*M)\otimes L_{\xi}}(V)
- {\rm ind}\, D_{F_{\xi},+}^{\wedge^{{\rm odd},0}(T^*M)\otimes L_{\xi}}(V)
\right) .\eqno (3.18)$$

{\bf Remark 3.12.}  Note that if $\xi$ is not contained in the image
   of the $\mu$, then it is automatically a regular value of $\mu$. In
   this case, $M_{G,\xi}=\emptyset$. By Theorem 0.2, one sees
   immediately that  in the summations in (3.17) and (3.18), the $\xi$ 
   actually runs through the integral lattice points contained in the image of $\mu$.

$\ $

Now if the  Euler characteristic $e(M)$ of $M$ is nonzero, 
then from (3.18) and
the Atiyah-Singer index theorem [AS], which gives that
$$\dim Q_{\rm dR}(M,L) =e(M),\eqno (3.19)$$
one deduces by (3.17) that at least one of the terms
$\dim Q(M, \wedge^{{\rm even},0}(T^*M)\otimes L) _{\xi} $ and 
$\dim Q(M, \wedge^{{\rm odd},0}(T^*M)\otimes L) _{\xi} $ 
should be nonzero.  In view of (3.16), (3.18), Theorem 3.6 and
Corollary 3.7, this  provides a concrete example mentioned in Remark 3.8.

$$\ $$

{\bf References}

$\ $

{\small 

[AS] M. F. Atiyah and I. M. Singer, The index of elliptic operators III. 
{\it Ann. of Math.} 87 (1968), 546-604.

[BL] J.-M. Bismut and G. Lebeau, Complex immersions and Quillen metrics. 
{\it Publ. Math. IHES.} 74 (1991), 1-297.

[DGMW] H. Duistermaart, V. Guillemin, E. Meinrenken and S. Wu,
Symplectic reduction and Riemann-Roch for circle actions.
{\it Math. Res. Lett.} 2 (1995), 259-266.

[G] V. Guillemin, Reduced phase space and Riemann-Roch. in
{\it Lie Groups and Geometry}, R. Brylinski et al. ed., Birkhaeuser
(Progress in Math., 123), 1995, pp. 305-334.

[GS] V. Guillemin and S. Sternberg, Geometric quantization and
multiplicities of group representations. 
{\it Invent. Math.} 67 (1982), 515-538.

[JK] L. C. Jeffrey and F. C. Kirwan, Localization and quantization
conjecture. {\it Topology} 36 (1997), 647-693.

[K] F. C. Kirwan, {\it Cohomology of Quotients in Symplectic and Algebraic
Geometry.}  Princeton Univ. Press, Princeton, 1984.

[Ko] B. Kostant, Quantization and unitary representations. {in} {\it Modern
Analysis and Applications.}  Lecture Notes in Math. vol. 170,
Springer-Verlag, (1970), pp. 87-207.

[M1] E. Meinrenken, On Riemann-Roch formulas for multiplicities.
{\it J. Amer. Math. Soc.} 9 (1996), 373-389.

[M2] E. Meinrenken, Symplectic surgery and the Spin$^c$-Dirac operator.
{\it  Adv. in Math.} To appear.

[MS] E. Meinrenken and R. Sjamaar, Singular reduction and quantization.
{\it Preprint}, 1996 (dg-ga 9707023).

[TZ1] Y. Tian and W. Zhang, An analytic proof of the geometric
quantization conjecture of Guillemin-Sternberg.
{\it Invent. Math.} To appear.

[TZ2] Y. Tian and W. Zhang, Quantization formula for symplectic manifolds
with boundary. {\it Preprint}, 1997.

[TZ3] Y. Tian and W. Zhang, Quantization formula for singular
reductions. {\it Preprint}, 1997 (dg-da 9706016).

[V1] M. Vergne, Multiplicity formula for geometric quantization, Part I.
{\it Duke Math. J}. 82 (1996), 143-179.

[V2] M. Vergne, Multiplicity formula for geometric quantization, Part II.
{\it Duke Math. J}. 82 (1996), 181-194.

$\ $

Youliang Tian, CUNY Graduate Center {\it and}  Courant Institute of Mathematical
               Sciences, New York, U.S.A.

{\it e-mail address:} ytian@cims.nyu.edu

$\ $

Weiping Zhang,
Nankai Institute of Mathematics,
Tianjin, 300071,
People's Republic of China

{\it e-mail address:} weiping@sun.nankai.edu.cn
\end{document}